\documentclass[twoside]{article}
\include{graphicx}
\title{The asymptotics of the Touchard polynomials: a uniform approximation}
\author{\sc R. B.\ Paris \\
{\em Division of Computing and Mathematics}, \\
{\em University of Abertay Dundee, Dundee DD1 1HG, UK}
}
\begin{document}
\def\f#1#2{\mbox{${\textstyle \frac{#1}{#2}}$}}
\def\dfrac#1#2{\displaystyle{\frac{#1}{#2}}}
\def\boldal{\mbox{\boldmath $\alpha$}}
{\newcommand{\Sgoth}{S\;\!\!\!\!\!/}
\newcommand{\bee}{\begin{equation}}
\newcommand{\ee}{\end{equation}}
\newcommand{\lam}{\lambda}
\newcommand{\ka}{\kappa}
\newcommand{\al}{\alpha}
\newcommand{\fr}{\frac{1}{2}}
\newcommand{\fs}{\f{1}{2}}
\newcommand{\g}{\Gamma}
\newcommand{\br}{\biggr}
\newcommand{\bl}{\biggl}
\newcommand{\ra}{\rightarrow}
\newcommand{\mbint}{\frac{1}{2\pi i}\int_{c-\infty i}^{c+\infty i}}
\newcommand{\mbcint}{\frac{1}{2\pi i}\int_C}
\newcommand{\mboint}{\frac{1}{2\pi i}\int_{-\infty i}^{\infty i}}
\newcommand{\gtwid}{\raisebox{-.8ex}{\mbox{$\stackrel{\textstyle >}{\sim}$}}}
\newcommand{\ltwid}{\raisebox{-.8ex}{\mbox{$\stackrel{\textstyle <}{\sim}$}}}
\renewcommand{\topfraction}{0.9}
\renewcommand{\bottomfraction}{0.9}
\renewcommand{\textfraction}{0.05}
\newcommand{\mcol}{\multicolumn}
\date{}
\maketitle
\pagestyle{myheadings}
\markboth{\hfill \sc R. B.\ Paris  \hfill}
{\hfill \sc  Touchard polynomials\hfill}
\begin{abstract}
The asymptotic expansion of the Touchard polynomials $T_n(z)$ (also known as the exponential polynomials)
for large $n$ and complex values of the variable $z$, where $|z|$ may be finite or allowed to be large like $O(n)$, has been recently considered in \cite{P1}. When $z=-x$ is negative, it is found that there is a coalesence of two contributory saddle points when $n/x=1/e$. Here we determine the expansion when $n$ and $x$ satisfy this condition and also a uniform two-term approximation involving the Airy function in the neighbourhood of this value.
Numerical results are given to illustrate the accuracy of the asymptotic approximations obtained.

\vspace{0.4cm}

\vspace{0.3cm}

\end{abstract}

\vspace{0.3cm}

\noindent $\,$\hrulefill $\,$

\vspace{0.2cm}

\begin{center}
{\bf 1. \  Introduction}
\end{center}
\setcounter{section}{1}
\setcounter{equation}{0}
\renewcommand{\theequation}{\arabic{section}.\arabic{equation}}
The Touchard polynomials $T_n(z)$, also known as exponential polynomials, are defined by
\bee\label{e11}
T_n(z)=e^{-z}\sum_{k=0}^\infty \frac{k^nz^k}{k!}=e^{-z} \bl(z \frac{d}{dz}\br)^n e^z
\ee
and were first introduced in a probabilistic context in 1939 by J. Touchard \cite{T}.
They have the generating function
\bee\label{e12}
\exp\,[z(e^t-1)]=\sum_{n=0}^\infty T_n(z)\,\frac{t^n}{n!}
\ee
and possess the alternative representation given by
\bee\label{e13}
T_n(z)=\sum_{k=0}^n S(n,k) z^k,
\ee
where $S(n,k)$ is the Stirling number of the second kind \cite[p.~624]{DLMF}.

In \cite{P1} we considered the asymptotic expansion of $T_n(z)$ for large $n$ and complex values of the variable $z$
by an application of the method of steepest descents applied to a contour integral representation. In this treatment $|z|$ was finite or allowed to be large like $O(n)$. It was found that there is an infinite number of saddle points of the integrand but that the precise number contributing to the expansion of $T_n(z)$ depended on the values of $n$ and $|z|$. When $z=-x$ ($x>0$), which is the central issue in this note, 
we have the expansions \cite[Theorem 2]{P1}
\bee\label{e14}
T_{n-1}(-x)\sim\left\{\begin{array}{ll}\Re \dfrac{\surd 2 \g(n)e^{x+n/t_0}}{\sqrt{\pi(1+t_0)}\ t_0^{n-1}} \sum_{s=0}^\infty \dfrac{c_{2s}(t_0) \g(s+\fs)}{n^{s+\fr} \g(\fs)} & (\mu>1/e)\\
\\
\dfrac{\g(n) e^{x+n/t_0}}{\sqrt{2\pi (1+t_0)}\,t_0^{n-1}}\sum_{s=0}^\infty\dfrac{c_{2s}(t_0) \g(s+\fs)}{n^{s+\fr} \g(\fs)} & (0<\mu<1/e) \end{array}\right.\ee
as $n\ra\infty$, where $t_0$ is one of the conjugate pair of roots of $te^t=-n/x$ with smallest modulus in the first expression and the smaller (negative) root in the second expression. Explicit expressions for the coefficients $c_{2s}(t_0)$ with $s\leq 2$ are given in \cite{P1}. In the case of the upper formula in (\ref{e14}), two conjugate saddles contribute to the expansion of $T_{n-1}(-x)$ when $1/e<\mu<\mu_1$, where $\mu_1\doteq 3.11179$; when $\mu\geq\mu_1$, there are other conjugate pairs of contributory saddles but these are not included in the upper formula in (\ref{e14}) as they as subdominant as $n\ra\infty$.

When $\mu:=n/x=1/e$, the two contributory saddle points coalesce to form a double saddle where the Poincar\'e-type expansions in (\ref{e14}) break down. In this note we obtain a uniform approximation for $T_{n-1}(-x)$ involving the Airy function together with an expansion valid when $\mu=1/e$.
Some numerical examples are given to illustrate the accuracy of the approximations obtained.

\vspace{0.6cm}

\begin{center}
{\bf 2. \ An integral representation}
\end{center}
\setcounter{section}{2}
\setcounter{equation}{0}
\renewcommand{\theequation}{\arabic{section}.\arabic{equation}}
From (\ref{e12}) we obtain the integral representation
\[
T_n(z)=\frac{n!\,e^{-z}}{2\pi i}\oint \frac{e^{ze^t}}{t^{n+1}}\,dt,
\]
where the integration path is a closed circuit described in the positive sense surrounding the origin. We let $z=-x$, where the variable $x>0$ is either finite or large like $O(n)$.
Since $|\exp (-xe^t)|\ra 0$ as $\Re (t)\ra+\infty$ when $|\Im (t)|<\fs\pi$, it follows that the closed path above may be opened up into a loop\footnote{In \cite{P1}, where $z$ is a complex variable and $n\geq 1$, the closed path around the origin was opened up into a loop which commences at $-\infty$, encirles the origin in the positive sense and returns to $-\infty$.}, which commences at $+\infty$, encircles the origin and returns to $+\infty$.
Then, introducing the scaled Touchard polynomial ${\hat T}_n(z)$ by
\[{\hat T}_n(z)\equiv \frac{1}{n!}T_n(z),\]
we have
\bee\label{e21}
{\hat T}_{n-1}(-x)=\frac{e^{x}}{2\pi i} \int_{\infty}^{(0+)} e^{n\psi(t)}dt,
\ee
where
\bee\label{e22}
\psi(t)\equiv\psi(t;\mu):=-\frac{e^{t}}{\mu}-\log\,t,\qquad 
\mu:=\frac{n}{x}~.
\ee 

Saddle points of the integrand occur when $\psi'(t)=0$; that is when
\[te^t=-\mu,\]
for which there is an infinite number of (complex) roots.  For a full discussion of the distribution of the saddle points see \cite{P1}.
When $0<\mu<1/e$, there are two saddles on the negative real axis given by the negative values of the Lambert-$W$ function;  see \cite[p.~111]{DLMF}. When $\mu=1/e$, these two saddles coalesce to form a double saddle point at $t=-1$ and when $\mu>1/e$ the saddles move off the real axis to form a complex conjugate pair.

In Fig.~1 we show examples of the steepest paths through the contributory saddles when (i) $0<\mu<1/e$, (ii) $1/e<\mu<\mu_1$ and (iii) $\mu=1/e$, where $\mu_1$ is specified in Section 1. The $t$-plane has a branch cut along $[0,\infty)$. 
In case (i), the saddles
%\footnote{In \cite{P1} these saddles were called $t_0$ and $t_1$, respectively.} 
$t_0$ and $t_1$ are situated on the negative real axis, with $t_0\in (0,-1)$ and $t_1\in (-1,-\infty)$ given by the negative branch of the Lambert-$W$ function \cite[p.~111]{DLMF}.
The paths of steepest descent emanating from $t_0$ pass to $+\infty$ and the paths of steepest ascent from $t_1$ asymptotically approach the lines $\Im (t)=\pm\pi$ as $\Re (t)\ra +\infty$. The integration path in (\ref{e21}) can then be deformed to pass over the steepest descent path emanating from $t_0$.
In case (ii), the saddles $t_0$ and $t_1$ form a conjugate pair and the integration path is the path labelled $ABCD$ in Fig.~1(b). When $\mu\geq\mu_1$, however, there are additional conjugate pairs of saddles (dependent on the value of $\mu$)
but these are subdominant as $n\ra\infty$.
In case (iii), the saddles $t_0$ and $t_1$ coalesce to form a double saddle at $t=-1$; the integration path then becomes the path $CSB$ in Fig.~1(c).
\begin{figure}[th]
	\begin{center}	{\tiny($a$)}\includegraphics[width=0.275\textwidth]{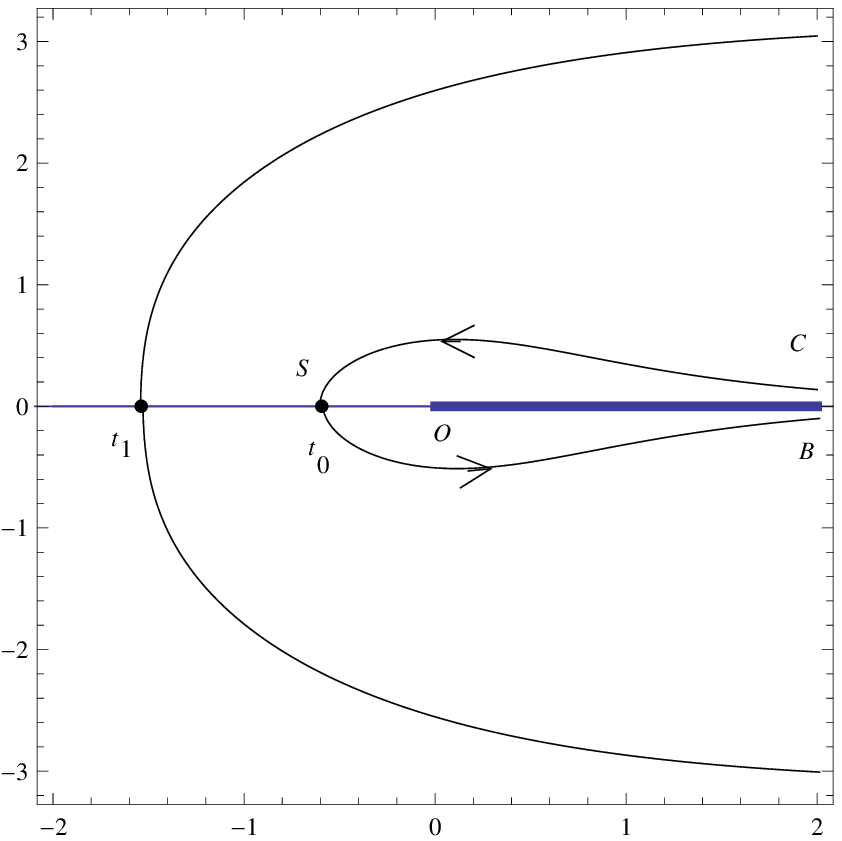}\qquad
	{\tiny($b$)}\includegraphics[width=0.275\textwidth]{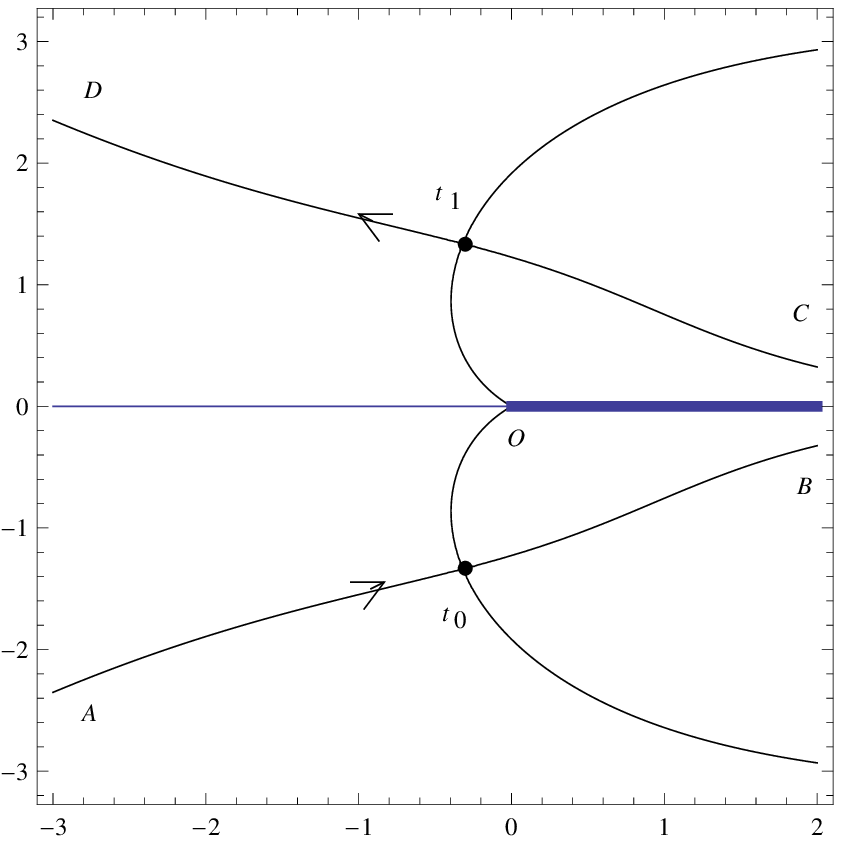}\qquad
	{\tiny($c$)}\includegraphics[width=0.275\textwidth]{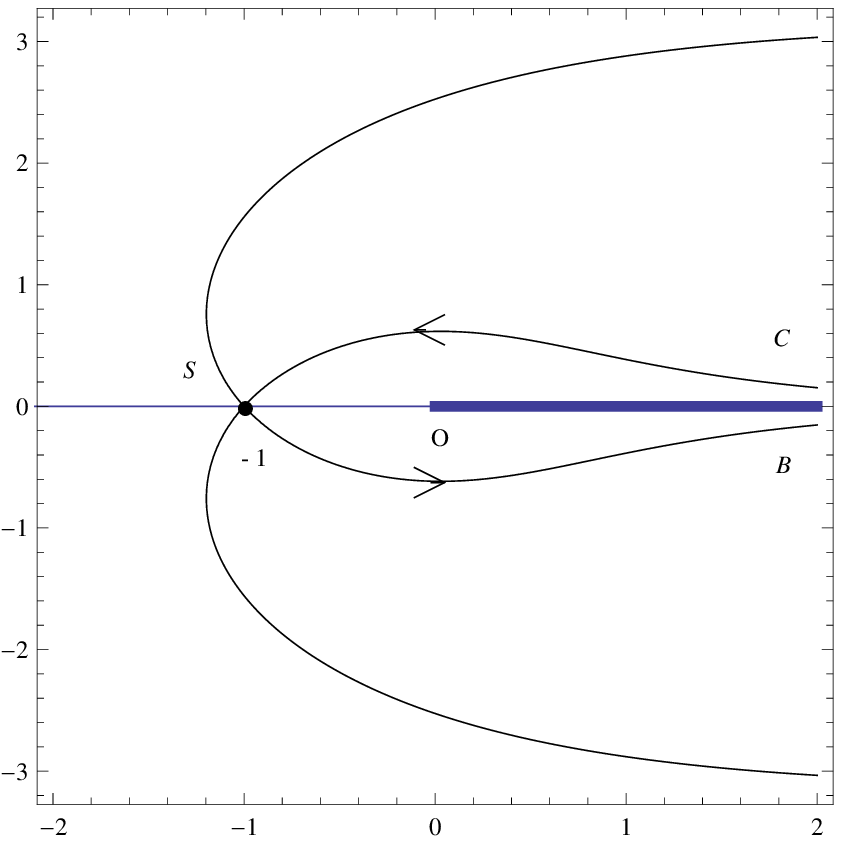}
\caption{\small{Paths of steepest descent and ascent through the saddles when (a) $0<\mu<1/e$, (b) $1/e<\mu<\mu_1$ and (c) $\mu=1/e$. The saddles are denoted by heavy dots; the arrows indicate the direction of integration  taken along steepest descent paths. There is a branch cut along $[0,\infty)$.}}
	\end{center}
\end{figure}
\vspace{0.6cm}

\begin{center}
{\bf 3. \ The asymptotics of ${\hat T}_{n-1}(-x)$ for $\mu\simeq 1/e$}
\end{center}
\setcounter{section}{3}
\setcounter{equation}{0}
\renewcommand{\theequation}{\arabic{section}.\arabic{equation}}
Both the expansions in (\ref{e14}) cease to be valid in the neighbourhood of the double saddle at $t=-1$.
We now determine an expansion valid at the coalescence point when $\mu=1/e$ and a uniform two-term approximation when $\mu \simeq 1/e$.
\vspace{0.3cm}

\noindent{\it 3.1\ \ The expansion of ${\hat T}_{n-1}(-x)$ when $\mu=1/e$}
\vspace{0.2cm}

\noindent
When $\mu=1/e$, the integration path can be deformed to coincide with the
steepest descent path ${\cal C}$ that enters $t=-1$ in the direction $\arg\,t=\pi/3$ and leaves to $t=-1$ in the direction $\arg\,t=-\pi/3$; see Fig.~1(c). If we put
\[-u=\psi(t)-\psi(-1)=\frac{\tau^3}{3!}+\frac{5\tau^4}{4!}+\frac{23\tau^5}{5!}+\frac{119\tau^6}{6!}+\frac{719\tau^7}{7!}+\cdots\ ,\qquad \tau=t+1\]
we find upon inversion using {\it Mathematica} that
\[\tau(w)=(6w)^{1/3}-\frac{5w^{2/3}}{2\cdot 6^{1/3}}+\frac{33w}{40}-\frac{1463w^{4/3}}{720\cdot 6^{2/3}}+\frac{126827w^{5/3}}{151200\cdot 6^{1/3}}-\frac{15451w^2}{44800}+ \ldots ,
\]
where $w=e^{-\pi i}u$ on the path $SB$ and $w=e^{\pi i}u$ on the path $SC$ in Fig.~1(c).
Then, upon differentiation of $\tau(w)$, we have
\begin{eqnarray*}
\frac{1}{2\pi i}\int_{\cal C} e^{-nu} \frac{d\tau}{du}\,du&=&\frac{1}{2\pi i} \int_0^\infty e^{-nu}\bl\{\frac{d\tau(ue^{-\pi i})}{du}-\frac{d\tau(ue^{\pi i})}{du}\br\}\,du\\
&=&\frac{1}{\pi}\int_0^\infty e^{-nu}\bl\{-\frac{6^{1/3} \sin \f{1}{3}\pi}{3 u^{2/3}}+\frac{5 \sin \f{2}{3}\pi}{3\cdot 6^{1/3} u^{1/3}}+\frac{1463 \sin \f{4}{3}\pi}{540\cdot 6^{2/3}} u^{1/3}+\cdots\br\} du  \\
&=&-\frac{1}{3\pi}\br\{\frac{\g(\f{1}{3}) \sin \f{1}{3}\pi}{(\f{1}{6}n)^{1/3}}-\frac{5 \g(\f{2}{3}) \sin \f{2}{3}\pi}{6 (\f{1}{6}n)^{2/3}}-\frac{1463 \g(\f{4}{3}) \sin \f{4}{3}\pi}{6480 (\f{1}{6}n)^{4/3}}- \cdots \br\}.
\end{eqnarray*}

Hence we obtain the following result.
\newtheorem{theorem}{Theorem}
\begin{theorem}$\!\!\!.$ Let $\mu=n/x=1/e$. Then as $n\ra\infty$, we have the expansion for the scaled Touchard polynomial
\bee\label{e38}
{\hat T}_{n-1}(-x)\sim (-)^{n-1} \frac{e^{x-n}}{3\pi} \sum_{m=0}^\infty \frac{(-)^m B_m \g(\f{1}{3}m+\f{1}{3})}{(\f{1}{6}n)^{(m+1)/3}}\,\sin \pi(\f{1}{3}m+\f{1}{3}),
\ee
where
\[B_0=1,\quad B_1=\frac{5}{6},\quad B_3=\frac{1463}{6480},\quad B_4=\frac{126827}{1088640},\quad B_6=\frac{4732223}{167961600}, \ldots\ .\]
\end{theorem}
We note the omission of the coefficients with index $m=2, 5, \ldots$; these terms do not contribute to the expansion on account of the vanishing of the sine factor.
\vspace{0.3cm}

\noindent{\it 3.2\ \ A uniform approximation for ${\hat T}_{n-1}(-x)$ when $\mu\simeq 1/e$}
\vspace{0.2cm}

\noindent
Let $\mu=1/(e\xi)$ with $\xi>0$; when $\xi\geq 1$ the saddles $t_0$ and $t_1$ are real, whereas when $\xi<1$ the saddles form a conjugate pair. To obtain a uniform approximation valid for $\xi\simeq 1$ we apply the standard cubic transformation \cite{CFU}
\bee\label{e21a}
\psi(t)=\f{1}{3}u^3-\zeta u+\beta
\ee
to the integrand in (\ref{e21}). The quantities $\zeta$ and $\beta$ depend on $\xi$ and are determined by the requirement that the saddles $t_0$ and $t_1$ correspond to $u=\zeta^{1/2}$ and $u=-\zeta^{1/2}$, respectively; that is
\bee\label{e21b}
\beta=\fs\{\psi(t_0)+\psi(t_1)\},\qquad \psi(t_j)=1/t_j-\log\,t_j\quad (j=0, 1)
\ee
and
\bee\label{e21c}
\f{2}{3}\zeta^{3/2}=\fs\{\psi(t_1)-\psi(t_0)\}\quad (\xi>1),\qquad \f{2}{3}(-\zeta)^{3/2}=\fs i\{\psi(t_0)-\psi(t_1)\}\quad (\xi<1).
\ee
For $x>0$, the quantity $\zeta\geq 0$ when $\xi\geq 1$ and $\zeta<0$ when $\xi<1$.

The integral for ${\hat T}_{n-1}(-x)$ in (\ref{e21}) then becomes
\[{\hat T}_{n-1}(-x)=\frac{e^{x+n\beta}}{2\pi i} \int_{\cal C'} e^{n(\frac{1}{3}u^3-\zeta u)} \frac{dt}{du}\,du,\]
where $\cal C'$ is the image in the $u$-plane of the integration path. With the substitution
\[g(u):=\frac{dt}{du}=A_0+B_0+(u^2-\zeta) G(u),\]
where
\bee\label{e23a}
A_0=\frac{1}{2}\{g(\zeta^\frac{1}{2})+g(-\zeta^\frac{1}{2})\},\qquad B_0=\frac{1}{2\zeta^\frac{1}{2}}\{g(\zeta^\frac{1}{2})-g(-\zeta^\frac{1}{2})\},
\ee
we then obtain \cite[p.~369]{W}, \cite[p.~67]{PH}
\begin{equation}\label{e24}
{\hat T}_{n-1}(-x)=e^{x+n\beta} \left\{\frac{A_0}{n^{1/3}} U(n^{2/3} \zeta)-\frac{B_0}{n^{2/3}} U' (n^{2/3}\zeta)+\frac{n^{-1}}{2\pi i}\int_{\cal C} e^{n(\frac{1}{3}u^3-\zeta u)} G'(u)\,du \right\},
\end{equation}
where the prime denotes differentiation with respect to the argument concerned and the function $U(z)$ is given by
\[U(z)=\frac{1}{2\pi i}\int_{\cal C'} e^{\frac{1}{3}\tau^3-z\tau}d\tau.\]
The path ${\cal C}'$ in the $u$-plane when $\xi>1$ and $\xi<1$ can be shown to be asymptotic to the rays $\arg\,u=\pm \pi/3$ traversed in the direction from the upper half-plane to the lower half-plane (we omit these details). From \cite[Eq.~(9.5.4)]{DLMF}, the function $U(z)$ is therefore given by the Airy function $-\mbox{Ai} (z)$.

From \cite[p.~367]{W}, we have
\[g(\zeta^\fr)=\bl(\frac{2\zeta^\fr}{\psi''(t_0)}\br)^{\!1/2},\qquad g(-\zeta^\fr)=\bl(\frac{-2\zeta^\fr}{\psi''(t_1)}\br)^{\!1/2}\qquad (\zeta\neq 0),\]
where $\psi''(t_j)=(1+t_j)/t_j^2$ ($j=0, 1$). Then, from (\ref{e23a}), we find 
\bee\label{e25a}
A_0=\frac{\zeta^\frac{1}{4}}{\sqrt{2}}\bl\{\bl(\frac{1}{\psi''(t_0)}\br)^{\!1/2}+\bl(\frac{-1}{\psi''(t_1)}\br)^{\!1/2}\br\}, \qquad B_0
=\frac{\zeta^{-\frac{1}{4}}}{\sqrt{2}}\bl\{\bl(\frac{1}{\psi''(t_0)}\br)^{\!1/2}-\bl(\frac{-1}{\psi''(t_1)}\br)^{\!1/2}\br\}
\ee
when $\xi>1$, and 
\bee\label{e25b}
A_0=\sqrt{2} |\zeta|^\frac{1}{4}\, \Re \bl[\bl(\frac{i}{\psi''(t_1)}\br)^{\!1/2}\br],\qquad
B_0=\sqrt{2} |\zeta|^{-\frac{1}{4}}\, \Im \bl[\bl(\frac{i}{\psi''(t_1)}\br)^{\!1/2}\br]
\ee
when $\xi<1$.
Hence, upon neglecting the third term in braces in (\ref{e24}) (which is o$(n^{-1})$), we obtain the following result.
\begin{theorem}$\!\!\!.$ Let $\mu=n/x=1/(e\xi)$, where $\xi>0$. Then we have the uniform two-term approximation for the scaled Touchard polynomial
\bee\label{e25}
{\hat T}_{n-1}(-x)\sim(-)^{n-1} e^{x+n\Re (\beta)}
\bl\{\frac{A_0}{n^{1/3}} {\mbox Ai} (n^{2/3} \zeta)-\frac{B_0}{n^{2/3}} {\mbox Ai}' (n^{2/3}\zeta)\br\}
\ee
as $n\ra\infty$. The quantities $\beta$ and $\zeta$ are defined in (\ref{e21b}) and (\ref{e21c}), where $\zeta\geq 0$ for $\xi\geq 1$ and $\zeta<0$ when $\xi<1$. The coefficients $A_0$ and $B_0$ are given in (\ref{e25a}) and (\ref{e25b}).
\end{theorem}

At coalescence when $\xi=1$ ($t=-1$, $u=0$) we have $A_0=g(0)=t'(0)$, $B_0=g'(0)=(t'(0))^2+t''(0)$, where $t(u)=dt/du$ and, by repeated differentiation of (\ref{e21a}),
\[t'(0)=\bl(\frac{2}{\psi'''(-1)}\br)^{\!1/3},\qquad t''(0)=-\frac{\psi^{iv}(-1)}{6\psi'''(-1)}\,\bl(\frac{2}{\psi'''(-1)}\br)^{\!2/3}.\]
Since $\psi'''(-1)=1$, $\psi^{iv}(-1)=5$, we obtain $A_0=2^{1/3}$, $B_0=-\f{5}{6}\cdot 2^{2/3}$. Use of the standard values
Ai($0)=3^{-2/3}/\g(\f{2}{3})$, Ai$\,'(0)=-3^{-1/3}/\g(\f{1}{3})$, together with $\psi(-1)=-1-\pi i$ (so that $\Re (\beta)=-1$), then shows that the approximation (\ref{e25}) at coalescence becomes
\bee\label{e28}
{\hat T}_{n-1}(-x)\sim (-)^{n-1} \frac{e^{x-n}}{3\pi}\bl\{\frac{\g(\f{1}{3}) \sin \f{1}{3}\pi}{(\f{1}{6}n)^{1/3}}+
\frac{5}{6}\,\frac{\g(\f{2}{3}) \sin \f{2}{3}\pi}{(\f{1}{6}n)^{2/3}}\br\}\qquad (\mu=1/e)
\ee
as $n\ra\infty$. This is seen to agree with the first two terms of the expansion in (\ref{e38}). 
\vspace{0.3cm}

%\newpage

\noindent{\it 3.3\ \ Numerical examples}
\vspace{0.2cm}

\noindent
In Table 1 we illustrate the accuracy of the expansion (\ref{e38}) by presenting values of the absolute relative error in ${\hat T}_{n-1}(-x)$ for different values of $n$ and truncation index $m$ at the coalescence point $\mu=1/e$. The value of ${\hat T}_{n-1}(-x)$ was computed from (\ref{e13}). Similarly, in Table 2, we show the absolute relative error
in ${\hat T}_{n-1}(-x)$ for different values of the coalescence parameter $\xi$ using the two-term approximation in (\ref{e25}) and, when $\xi=1$, using (\ref{e28}).
\begin{table}[th]
\caption{\footnotesize{Values of absolute relative error in the computation of ${\hat T}_{n-1}(-x)$ using the asymptotic expansion (\ref{e38}) for different $n$ and truncation index $m$ when $\mu=1/e$.}}
\begin{center}
\begin{tabular}{l|lll}
\hline
&&&\\[-0.25cm]
\mcol{1}{c|}{$m$} & \mcol{1}{c}{$n=50$} & \mcol{1}{c}{$n=80$} & \mcol{1}{c}{$n=121$}\\
[0.1cm]\hline
&&&\\[-0.25cm]
0 & $2.514\times 10^{-1}$ & $2.095\times 10^{-1}$ & $1.788\times 10^{-1}$\\
1 & $8.558\times 10^{-3}$ & $5.390\times 10^{-3}$ & $3.585\times 10^{-3}$\\
3 & $2.744\times 10^{-3}$ & $1.437\times 10^{-3}$ & $8.144\times 10^{-4}$\\
4 & $1.638\times 10^{-4}$ & $6.490\times 10^{-5}$ & $2.868\times 10^{-5}$\\
6 & $6.184\times 10^{-5}$ & $2.029\times 10^{-5}$ & $7.616\times 10^{-6}$\\
[0.1cm]\hline
\end{tabular}
\end{center}
\end{table}
\begin{table}[th]
\caption{\footnotesize{Values of absolute relative error in the computation of ${\hat T}_{n-1}(-x)$ using the uniform approximation (\ref{e25}) and (\ref{e28}) for different values of the coalescence parameter $\xi$.}}
\begin{center}
\begin{tabular}{l|ll||l|ll}
\hline
&&&&&\\[-0.25cm]
%\mcol{1}{c|}{} & \mcol{2}{c|}{$n=20$} & \mcol{2}{c|}{$n=30$} & \mcol{2}{c}{$n=50$}\\
\mcol{1}{c|}{$\xi$} & \mcol{1}{c}{$n=81$} & \mcol{1}{c||}{$n=100$} & \mcol{1}{c|}{$\xi$} & \mcol{1}{c}{$n=81$} & \mcol{1}{c}{$n=100$} \\
[0.1cm]\hline
&&&&&\\[-0.25cm]
0.80 & $5.243\times 10^{-3}$ & $8.179\times 10^{-3}$ & 1.01 & $5.300\times 10^{-3}$ & $4.301\times 10^{-3}$\\
0.90 & $7.413\times 10^{-3}$ & $3.322\times 10^{-3}$ & 1.05 & $5.204\times 10^{-3}$ & $4.222\times 10^{-3}$\\
0.95 & $5.545\times 10^{-3}$ & $4.540\times 10^{-3}$ & 1.10 & $5.122\times 10^{-3}$ & $4.153\times 10^{-3}$\\
0.99 & $5.356\times 10^{-3}$ & $4.355\times 10^{-3}$ & 1.20 & $5.010\times 10^{-3}$ & $4.060\times 10^{-3}$\\
1.00 & $5.324\times 10^{-3}$ & $4.326\times 10^{-3}$ & 1.40 & $4.878\times 10^{-3}$ & $3.951\times 10^{-3}$\\
[.1cm]\hline
\end{tabular}
\end{center}
\end{table}
\vspace{0.6cm}

%\noindent{\bf Acknowledgement:}\ \ The author wishes to acknowledge V. Vinogradov for bringing the problem of the asymptotics of the Touchard polynomials to his attention.

\end{document}